\documentclass{amsart}

\begin{document}

\title[Canonical lifts]{The  Euler top and canonical lifts}
\bigskip

\def \h{\hat{\ }}
\def \cO{\mathcal O}
\def \ra{\rightarrow}
\def \bZ{{\mathbb Z}}
\def \cP{\mathcal V}
\def \cH{{\mathcal H}}
\def \cB{{\mathcal B}}
\def \d{\delta}
\def \cC{{\mathcal C}}
\def \jor{\text{jor}}

\newtheorem{THM}{{\!}}[section]
\newtheorem{THMX}{{\!}}
\renewcommand{\theTHMX}{}
\newtheorem{theorem}{Theorem}[section]
\newtheorem{corollary}[theorem]{Corollary}
\newtheorem{lemma}[theorem]{Lemma}
\newtheorem{proposition}[theorem]{Proposition}
\newtheorem{thm}[theorem]{Theorem}
\theoremstyle{definition}
\newtheorem{definition}[theorem]{Definition}
\theoremstyle{remark}
\newtheorem{remark}[theorem]{Remark}
\newtheorem{example}[theorem]{\bf Example}
\numberwithin{equation}{section}
\address{Department of Mathematics and Statistics\\University of New Mexico \\
  Albuquerque, NM 87131, USA\\ 
   Department of Mathematics and
  Statistics\\Boston University \\Boston, MA 02215, USA}
\email{buium@math.unm.edu \\ep@bu.edu } 
\maketitle

\bigskip

\medskip
\centerline{\bf Alexandru Buium and Emma Previato}
\bigskip

\begin{abstract} In this note, we prove 
a finiteness result for fibers that are canonical lifts in a given elliptic
fibration. 
The question
  was motivated by the authors' construction of an arithmetic Euler top,
  and it highlights an interesting discrepancy between the arithmetic and the
  classical case: in the former, it is impossible to extend the flows to a
  compactification of the phase space, viewed as an elliptic fibration over
  the space of action variables.
\end{abstract}

\section{Introduction}

In this note, we present an analogue, for the multiplicative group ${\mathbb G}_m$, of one of the main finiteness results for elliptic curves proved 
in \cite{local}. This is then applied to 
 a question left open in \cite[Remark 6.6]{euler}, namely
 whether 
the main theorem, \cite[Th. 6.1]{euler}
could be improved by showing that a certain extension property 
(saying, roughly, that an {\it arithmetic flow} on an affine elliptic fibration can be ``compactified") holds
for the {\it arithmetic  Euler top}, which was there introduced.
 As we shall see in Corollary \ref{tootoo} below the answer to this question is negative: no such extension property holds and, hence, no such improvement of Theorem 6.1 in \cite{euler} is possible.  
Note, on the other hand,   that  
the corresponding extension property {\it does hold} in the case of the {\it
  classical} Euler top,
which is an algebraically completetly integrable Hamiltonian system (ACI).
 This discrepancy between the arithmetic and the
classical case adds a somewhat intriguing feature to the
construction effected in \cite{euler}.   

Although the  motivation of the present paper arose from \cite{euler}, 
 the  paper is written so as to be logically independent.

In section 2 we discuss terminology and conventions and we state two
theorems. The first, Theorem \ref{thm1}, is the aforementioned
 analogue, for the multiplicative group ${\mathbb
  G}_m$, of a finiteness result for elliptic curves proved  
in \cite{local}. The second, Theorem \ref{thm2},
is an application of the first theorem to
the geometric setting of the Euler top. We then derive a corollary to 
Theorem \ref{thm2},
 showing the relevance of our analysis to the  question  raised
in \cite{euler}.  Section 3 contains the proofs of the two theorems and a
final remark on a variation  on Theorem \ref{thm2}. 

\medskip

{\bf Acknowledgments}. The authors are indebted to the IHES in Bures sur Yvette, where their collaboration on this project started.
The first author acknowledges partial support from the Max Planck Institute for Mathematics in Bonn, from  the Simons Foundation (award 311773), and from
NSF grant DMS-1606334 that funded the DART VII conference.
The second author expresses her sincerest
 gratitude to the Boston University
College of Arts and Sciences Associate Dean Stan Sclaroff and Assistant Dean
Richard Wright, for the travel support that made this research possible.

\section{Main results}

First, some terminology. In what follows varieties over ${\mathbb C}$ will be
identified with their sets 
of ${\mathbb C}$-points. By an {\it elliptic fibration} over  ${\mathbb C}$ we
understand a surjective morphism $f:{\mathcal E}\ra B$ from a non-singular
projective surface to a non-singular projective curve, whose general fiber is
a connected curve of genus $1$. If $B^0\subset B$ is the set of all points
with smooth fiber then one has the $J$-function $J:B^0\ra {\mathbb C}$ 
that attaches to each $P\in B^0$ the $j$-invariant $j({\mathcal E}_P)$ of the
fiber ${\mathcal E}_P=f^{-1}(P)$.  
The theory of $J$ is usually presented for {\it Jacobian fibrations}, i.e.,
for elliptic fibrations possessing a section. On the other hand, 
to any elliptic fibration  $f:{\mathcal E}\ra B$ one can attach its {\it
  Jacobian fibration} $f':{\mathcal E}'\ra B$ which is a compactification of 
$Pic^0_{{\mathcal E}^0/B^0}\ra B^0$ 
where ${\mathcal E}^0:=f^{-1}(B^0)$; cf. \cite{neron}, p. 260. Then the fibers
of $f$ and $f'$ above each $P\in B^0$ are isomorphic so $f$ and $f'$ have the
same $J$-function $B^0\ra {\mathbb C}$. This $J$-function induces a morphism 
 $J:B\ra {\mathbb P}^1={\mathbb C}\cup \{\infty\}$ which we still refer to as
the {\it $J$-function}, cf. \cite{bogo}.

 We also make, in what follows, the following convention. If finitely many
 varieties and morphisms over ${\mathbb C}$ are given, which are defined over
 $\overline{\mathbb Q}$, then we will implicitly choose a number field $L$
 over which they are defined and we will choose models of our objects over
 $\cO_L[1/M]$ for some integer $M$. Then, for any sufficiently big rational
 prime $p>>0$, we will  
 choose an unramified prime ${\mathfrak P}$ in $\cO_L$, not containing $M$,
 and we will 
 base change these models to $R=R_{\mathfrak P}$, the completion of the
 maximum unramified extension of the ${\mathfrak P}$-adic  completion of
 $\cO_L$, which is isomorphic to the completion $R_p$
 of the maximum unramified  
extension of $\bZ_p$. These new schemes and morphisms over $R_p$ will be
 denoted by the same letters as our original varieties and morphisms over
 ${\mathbb C}$. We will denote by $k_p=R_p/pR_p$ the residue field
 of $R_p=R_{\mathfrak P}$, which is therefore
 an algebraic closure of the field
 ${\mathbb F}_p=\bZ/p\bZ$. Also we denote by $\ \widehat{\ }\ $ 
the $p$-adic
 completion of  rings or schemes. 

Further, we denote by ${\mathbb G}_m=Spec\ \bZ[t,t^{-1}]$ the multiplicative group scheme over the integers; so for any ring $A$ we have ${\mathbb G}_m(A)=A^{\times}$.
According to our conventions,
over ${\mathbb C}$, we continue to write ${\mathbb G}_m$
in place of ${\mathbb G}_m({\mathbb C})={\mathbb C}^{\times}$ and 
 we have ${\mathbb C}={\mathbb G}_m\cup\{0\}$ and ${\mathbb
  P}^1={\mathbb G}_m\cup \{0,\infty\}$. 

Finally recall that 
a Frobenius lift on a scheme (respectively on a $p$-adic formal scheme) is an
endomorphism  whose reduction mod $p$ is the absolute Frobenius.  
Trivially, $Spec\ R_p$ has a unique Frobenius lift.
An elliptic curve over $R_p$ is said to be
a {\it canonical lift} (or to be CL) if it has a Frobenius lift which is
compatible with the Frobenius lift on $Spec\ R_p$. 
If an elliptic curve is CL then it has complex
multiplication (or is CM); 
cf. \cite{local} for a review of the relationship  between CL and CM.

Our first result will be the following:

\begin{theorem}
\label{thm1}
Let $f:{\mathcal E}\ra {\mathbb P}^1={\mathbb G}_m\cup\{0,\infty\}$ be an
elliptic fibration defined over $\overline{\mathbb Q}$, with $J$-function
$J:{\mathbb P}^1\ra {\mathbb P}^1$. Assume there exists $P_0\in {\mathbb G}_m$
such that  $J(P_0)=\infty$. Then for all sufficiently big primes $p>>0$ 
there exists a constant $c_p$  such that for any integer $r$ and any subgroup
$\Gamma\subset {\mathbb G}_m(R_p)=R_p^{\times}$ of  rank $r<\infty$ the set 
$$\{P\in \Gamma;\ {\mathcal E}_P\ \ \text{is}\ \ CL\}$$
is finite of cardinality at most $c_p p^r$. 
\end{theorem}

Here the rank of $\Gamma$  is defined, as usual,
 as the maximum number of
multiplicatively independent elements in $\Gamma$; in particular, the torsion
group $T_p$ of $R_p^{\times}$ has rank $0$. Also, according to our
conventions, ${\mathcal E}_P$ in the above theorem is the fiber $f^{-1}(P)$ 
of $P\in {\mathbb P}^1(R_p)$,
 where $f:{\mathcal E}\ra {\mathbb P}^1$ is the induced morphism over $R_p$. 

What we will need in the application of Theorem \ref{thm1} to the next result
will be the finiteness statement, not an actual estimate for the
cardinality.

\medskip

To state our next result recall from \cite{euler} the geometric setting of the
Euler equations. Let $a_1,a_2,a_3\in {\mathbb C}$ be distinct complex numbers
and consider the functions
\begin{equation} 
\label{HH}
H_1=\sum_{i=1}^3a_ix_i^2,\ \ \ H_2=\sum_{i=1}^3x_i^2,\end{equation}
in the polynomial ring ${\mathbb C}[x_1,x_2,x_3]$.
Also we consider the affine spaces 
\begin{equation}
\label{condition 3}
{\mathbb A}^2=\operatorname{Spec}\ {\mathbb C}[z_1,z_2],\ \ \ {\mathbb A}^3=\operatorname{Spec}\ {\mathbb C}[x_1,x_2,x_3]
\end{equation}
and the morphism
\begin{equation}
\label{condition 4}
{\mathcal H}:{\mathbb A}^3\ra {\mathbb A}^2
\end{equation}
defined by
\begin{equation}
z_1\mapsto H_1,\ \ \ z_2\mapsto H_2.\end{equation}
 For any ${\mathbb C}$-point $c=(c_1,c_2)\in {\mathbb C}^2={\mathbb
 A}^2({\mathbb C})$ we denote by  
\begin{equation}
\label{definition of Ec}
E_c=\operatorname{Spec}\ {\mathbb C}[x_1,x_2,x_3]/(H_1-c_1,H_2-c_2)
\end{equation}
the fiber of ${\mathcal H}$ at $c$.
Consider the polynomial
\begin{equation}
 \label{definition of N}
 N(z_1,z_2)= \prod_{i=1}^3(z_1-a_iz_2)\in {\mathbb C}[z_1,z_2].\end{equation}
 For 
 $c\in {\mathbb A}^2$ with $N(c)\neq 0$,
 $E_c$ is smooth over ${\mathbb C}$. 
If we consider  the  projective closure ${\mathcal E}_c$ of $E_c$ in the
 projective space $${\mathbb P}^3=Proj\ {\mathbb C}[t_0,t_1,t_2,t_3],\ \ \
 x_i=t_i/t_0,$$ 
then for $c\in {\mathbb A}^2$ with 
$N(c)\neq 0$
the curve ${\mathcal E}_c$ is still smooth.
(In \cite{euler}, the curves ${\mathcal E}_c$ were denoted by $E^*_c$.)

\medskip

Recall that for a rational prime $p$ we denoted by $T_p$ the torsion subgroup
of $R_p^{\times}$; it is the set of all roots of unity of order prime to $p$
in an algebraically closed field containing $R_p$.  
In the notation above we will prove the following:

\begin{theorem}\label{thm2}
Let $a_1,a_2,a_3\in \overline{\mathbb Q}$ and $c_2\in{\mathbb C}$ a root
of unity. For all sufficiently big primes $p>>0$  the set 
$$\{c_1\in T_p;\ \ N(c_1,c_2)\in R_p^{\times}, \ {\mathcal E}_{(c_1,c_2)}\ \
\text{is}\ \ CL\}$$ 
is finite.
\end{theorem}

To explain the relevance of Theorem \ref{thm2} for the  Euler equations
 considered in \cite{euler}, we  first recall the classical picture. 
 The classical {\it Euler top} 
is a rotating solid body attached to a fixed point in 
(three-dimensional) space,
subject to no external force. Its motion is described by   a system of three ordinary
 (non-linear) differential equations in $3$ variables (the {\it Euler equations}); these equations correspond to a
 polynomial vector field on the affine complex $3$-space and hence to a
 derivation $\d$ on the polynomial ring ${\mathbb C}[x_1,x_2,x_3]$. The
 derivation is given by the expression 
\begin{equation}
\label{classical EF}
\d:=(a_2-a_3)x_2x_3\frac{\partial}{\partial x_1}+(a_3-a_1)x_3x_1\frac{\partial}{\partial x_2}+(a_1-a_2)x_1x_2\frac{\partial}{\partial x_3},\end{equation}
where $a_1,a_2,a_3\in {\mathbb C}$ are distinct complex numbers. This vector field is trivially seen to have
 $H_1$ and $H_2$ in \ref{HH} as 
  prime integrals 
in the sense that
\begin{equation}
\label{HH0}
 \d H_1=\d H_2=0.\end{equation}
For  $c=(c_1,c_2)\in {\mathbb C}^2$ with $N(c)\neq 0$ the
curves $E_c$  given 
by \ref{definition of Ec} are  tangent to the vector field so  the derivation
$\d$ induces derivations $\d_c$ on the rings $\cO(E_c)$. The remarkable fact
about the situation is that  $\d_c$ on $\cO(E_c)$ are   {\it linearized} in
the sense   that for all $c\in {\mathbb A}^2$ with $N(c)\neq 0$, 
\begin{equation}
\label{fomist}
\textit{$\d_c$ extends to a vector field on the compactification ${\mathcal
    E}_c$ of $E_c$.}\end{equation} 
This makes the Euler equations  ``solvable by elliptic functions."
Condition \ref{fomist} is equivalent  to the condition that
 \begin{equation}
  \label{Lie}
  \d_c\omega_c=0,\end{equation}
  where we continue to  denote by $\d_c$ the action of $\d_c$ as Lie
  derivative on the $1$-forms on $E_c$ and  $\omega_c$ is the restriction to
  $E_c$ of some (equivalently any) non-zero $1$-form on ${\mathcal E}_c$; the equivalence of
  \ref{fomist} and \ref{Lie} is, of course, a consequence of Henri Cartan's   formula for the Lie derivative.
  
  \medskip
   
  In \cite{euler}, we developed an arithmetic analogue of the classical
  Euler equations. 
  To explain this, it is convenient to introduce some ad hoc terminology. 
  Let us define a
    {\it $p$-triple} as being a triple $(K,X,\phi)$ where
   
   \medskip
  
  $\bullet$ $K\in \cO(\widehat{{\mathbb A}^2})=\widehat{R_p[z_1,z_2]}$,  $K\not\equiv 0$ mod $p$,
    
  $\bullet$  
  $X\subset {\mathbb A}^3$ an open set over $R_p$, 
  
  $\bullet$ $\phi$   a Frobenius lift  on  $\widehat{X}$.
  
  \medskip

  \noindent such that the reduction mod $p$ of $K$ is a homogeneous
  polynomial in $k_p[z_1,z_2]$,

\medskip
  
 Theorem 6.1  in \cite{euler} implies the following:
   
   \begin{corollary}
   \label{rory}
  Let $a_1,a_2,a_3\in \overline{\mathbb Q}$. Then, if $p>>0$ is any sufficiently big prime, there exists  a $p$-triple $(K,X,\phi)$ satisfying the following two conditions:
  
  1) one has equalities
  \begin{equation}
\label{p1}
\phi(H_1)=H_1^p,\ \ \phi(H_2)=H_2^p;\end{equation}

2) for all $c\in T_p^2$ with $N(c)K(c)\in R_p^{\times}$, one has
 $\widehat{E}_c\cap \widehat{X}\neq \emptyset$ and
 \begin{equation}
 \label{p2}
 \frac{\phi_c^*}{p}\omega_c\equiv K(c)^{-1}\cdot \omega_c\ \ \ \text{mod}\ \ \ p,\end{equation}
    where 
   
   \medskip
   
      $\phi_c$ is  the Frobenius lift on  $\widehat{E}_c\cap \widehat{X}$ induced by $\phi$;
   
     $\omega_c$ is an $R_p$-basis for the $1$-forms on the compactification ${\mathcal E}_c$ of $E_c$ over $R_p$.\end{corollary}
 
 Here $\phi$ induces a Frobenius lift on $\widehat{E}_c\cap \widehat{X}$ because, for 
 $c=(c_1,c_2)\in T_p^2$, we have
 $$\phi(H_i-c_i)=H_i^p-c_i^p\in (H_i-c_i),\ \ \ i=1,2.$$
   We proved Corollary 
   \ref{rory} in \cite{euler} by choosing $K$ in the $p$-triple 
    to be the {\it Hasse invariant} $A_{p-1}$ of an appropriate associated
   plane quartic, cf. \cite{euler}; note that $A_{p-1}$ itself is a
   homogeneous polynomial of degree $p-1$ and 
   for $c_1,c_2\in \bZ_p$ we have that $A_{p-1}(c_1,c_2)\in \bZ_p^{\times}$ if
   and only if ${\mathcal E}_{(c_1,c_2)}$ has ordinary reduction. 
   
    Conditions \ref{p1} and \ref{p2} above are of course to be viewed as
    arithmetic analogues of conditions \ref{HH0} and \ref{Lie} respectively.  
    Now,  in view of the equivalence between \ref{fomist} and \ref{Lie}, it is
    natural to ask for an arithmetic analogue of the condition
    \ref{fomist}. An arithmetic analogue of \ref{fomist} could be the
    condition that    
      \begin{equation}
    \label{moore}
    \textit{$\phi_c$ extends to an endomorphism of 
    the compactification ${\mathcal E}_c$ of $E_c$.}
    \end{equation}
    Note that condition \ref{moore} implies, of course, the condition that
    \begin{equation}
    \label{remi}
    \frac{\phi_c^*}{p}\omega_p=\kappa_c\cdot \omega_c
    \end{equation}
    for some $\kappa_c\in R_p$ where $\omega_c$ is a basis for the space of
    $1$-forms on ${\mathcal E}_c$; in its turn condition \ref{remi} is a
    strengthening of  
    the congruence \ref{p2} (at least if the ``eigenvalue" $K(c)^{-1}$ in
    \ref{p2} is not being specified).  

    One is then tempted to ask the following question; cf.
   Remark 6.6 of \cite{euler}:
     \begin{equation}
    \textit{\it Does Corollary \ref{rory} 
    hold with condition \ref{p2} replaced by condition
    \ref{moore}?}\end{equation}
    The answer to this question is {\it no}; 
    indeed we have the following consequence of Theorem \ref{thm2}:

    \begin{corollary}
    \label{tootoo}
    Let $a_1,a_2,a_3\in \overline{\mathbb Q}$. Then, if $p>>0$ is any
    sufficiently big prime, there is no $p$-triple $(K,X,\phi)$  satisfying
    the following two conditions: 
    
    1) one has equalities
  \begin{equation}
\label{p11}
\phi(H_1)=H_1^p,\ \ \phi(H_2)=H_2^p;\end{equation}

 2) for all $c\in T_p^2$ with $N(c)K(c)\in R_p^{\times}$, one has
 $\widehat{E}_c\cap \widehat{X}\neq \emptyset$ and
 \begin{equation}
 \label{p22}
  \textit{$\phi_c$ extends to an endomorphism of 
    the compactification ${\mathcal E}_c$ of $E_c$,}
   \end{equation}
  where   $\phi_c$ is  the Frobenius lift on  $\widehat{E}_c\cap \widehat{X}$ induced by $\phi$.
    \end{corollary}

    {\it Proof}. Assume that there is an infinite set $S$ of primes  $p$ such
        that for each $p\in S$ there is a $p$-triple $(K,X,\phi)$  satisfying
        \ref{p11} and \ref{p22}. Let us fix a root of unity $c_2\in {\mathbb
        C}$. 
        By Theorem \ref{thm2} for all $p>>0$   the set 
   $$\{c_1\in T_p;\ \ N(c_1,c_2)\in R_p^{\times}, 
   \ {\mathcal E}_{(c_1,c_2)}\ \ \text{is}\ \ CL\}$$
   is finite.  Now,  for any $p\in S$, since the image of $K$ 
   in $k_p[z_1,z_2]$
   is homogeneous and not $\equiv 0$ mod $p$,
   the zero locus of $K$ mod $p$ in the plane ${\mathbb A}^2(k_p)=k_p^2$ is a
   union of lines passing through the origin hence the intersection of this
   locus with the line $z_2=c_2$ must be finite. So, for all $p\in S$, 
   the set
   $$\{c_1\in T_p;\ K(c_1,c_2)\not\in R_p^{\times}\}$$
   is finite. Similarly, since $N$ is homogeneous, the set
   $$\{c_1\in T_p;\ N(c_1,c_2)\not\in R_p^{\times}\}$$
   is finite for all $p>>0$. 
    It follows that for all except
   finitely many $p\in S$ the set 
   $$\{c_1\in T_p;\ \ N(c_1,c_2)K(c_1,c_2)\in R_p^{\times}, 
   \ {\mathcal E}_{(c_1,c_2)}\ \ \text{is not}\ \ CL\}$$
   has a finite complement in $T_p$, in particular it is non-empty. This
   violates  
   the condition 2) in the statement of the corollary. We obtained a
   contradiction which ends our proof. \qed 
   
   \begin{remark}
   As we saw, Corollary \ref{rory} fails if one replaces congruence \ref{p2} by condition \ref{moore}. One can ask, however, 
   whether  Corollary \ref{rory} continues to hold if one replaces  congruence \ref{p2} by the equality,
   \begin{equation}
   \label{p3}
   \frac{\phi_c^*}{p}\omega_c= K(c)^{-1}\cdot \omega_c,
   \end{equation}
   or, at least, by condition \ref{remi}.
   Indeed, \ref{p3} (an equality)
is stronger than \ref{p2} (the corresponding congruence),
and also stronger than condition \ref{remi}; on the other hand, as already
   noticed, condition \ref{remi} is weaker than condition  
   \ref{moore}. 
   \end{remark}

\section{Proof of the Theorems}

{\it Proof of Theorem \ref{thm1}}. 
Consider the 
restriction $J:{\mathbb G}_m\ra X(1)={\mathbb P}^1$ of the map $J:{\mathbb P}^1\ra {\mathbb P}^1$ in the statement of the Theorem and consider the
cartesian diagram
$$
\begin{array}{rcl}
{\mathbb G}_m & \stackrel{J}{\longrightarrow} & X(1)\\
\uparrow & \  & \uparrow\\
Z & \longrightarrow & X_1(N)
\end{array}
$$
where we are using the standard notation $X(1)$ and $X_1(N)$ for the modular curves, $N\geq 4$;
cf. \cite{local}. Note that since the right arrow is finite we have that the left arrow is finite so $Z$ is an affine curve. Recall that $X_1(N)$ possesses a distinguished cusp $\infty$ lying over the point $\infty\in X(1)={\mathbb P}^1$. So $Z$ contains a point $z_{\infty}$ mapping to $P_0\in {\mathbb G}_m$ and also to $\infty\in X_1(N)$. Denote by $X$ the normalization of an irreducible component of $Z_{red}$ that contains $z_{\infty}$, let $x_{\infty}\in X$ be a point lying over $z_{\infty}$, and consider the induced maps
$$\Pi:X\ra Z\ra X_1(N)$$
and 
$$\Phi:X\ra Z\ra {\mathbb G}_m\stackrel{L_{P_0}^{-1}}{\longrightarrow} {\mathbb G}_m,$$
where $L_{P_0}$ is the translation by $P_0$ in the group ${\mathbb G}_m$.
So we have at our disposal a diagram (a correspondence)
\begin{equation}
\label{latelate}
X_1(N) \stackrel{\Pi}{\longleftarrow} 
X \stackrel{\Phi}{\longrightarrow} {\mathbb G}_m,\ \ \ 
\Pi(x_{\infty})=\infty,\ \ \ \Phi(x_{\infty})=1.
\end{equation}
 Also $\Phi$ is a finite morphism. Note that for $p>>0$ we have
that any point $P$ in
$$
\{P\in \Gamma;\ {\mathcal E}_P\ \ \text{is}\ \ CL\}
$$
which is not a ramification point for $\Phi$ is the image of a point in
\begin{equation}
\label{int}
\Phi^{-1}(\Gamma')\cap \Pi^{-1}(CL)\end{equation}
where $\Gamma':=\langle \Gamma,P_0\rangle$ is the subgroup of ${\mathbb G}_m(R_p)=R_p^{\times}$ generated by $\Gamma$ and $P_0$, while $CL\subset X_1(N)(R_p)$ is the set of $CL$ points (points corresponding to CL elliptic curves).
Since $rank(\Gamma')\leq rank(\Gamma)+1$ it is sufficient to prove that the intersection \ref{int} has cardinality bounded by a constant that does not depend on $\Gamma'$ times $p^{rank(\Gamma')}$. 
This is a statement analogous to that of Theorem 1.6 in \cite{local}
and can be given a proof entirely analogous to the proof in \cite{local}. \qed
 
 \begin{remark}
 Some comments 
 on the proof of
Theorem 1.6 of \cite{local} and its relation to our situation here are in order. That theorem has two cases. 

One case refers to correspondences 
 \begin{equation}
\label{aaron}
X_1(N) \stackrel{\Pi}{\longleftarrow} 
X \stackrel{\Phi}{\longrightarrow} A,\ \ \ 
\Pi(x_{\infty})=\infty,\ \ \ \Phi(x_{\infty})=0,
\end{equation}
with $A$ an elliptic curve. The proof  in \cite{local} for this case is based on Fourier expansions and is divided into two subcases: one subcase corresponds to the situation when  $A$ itself is  CL while another subcase corresponds to the situation when  $A$ is not CL. The proof of the subcase in which $A$ is CL 
  goes through, essentially word for word,  when one replaces \ref{aaron} by  
\ref{latelate} ; all one has to do is  replace the canonical $\d$-character  of $A$ in section 4.3 of \cite{local}, cf. also \cite{char},  with the canonical $\d$-character of ${\mathbb G}_m$ in \cite{char}.

 Another case of Theorem 1.6 in \cite{local}  refers to correspondences 
  \begin{equation}
\label{land}
S \stackrel{\Pi}{\longleftarrow} 
X \stackrel{\Phi}{\longrightarrow} A,\ \ \ 
\ \ \ \Phi(x_{\infty})=0,
\end{equation}
with $A$ an elliptic curve and $S$ a Shimura curve.
 The proof in \cite{local}, in this case, uses Serre-Tate expansions (rather than Fourier expansions) and, interestingly, the Serre-Tate expansion method does not seem to go through for correspondences \ref{latelate}. Indeed 
 the Serre-Tate expansion method in \cite{local} only applied 
 to primes that were not {\it anomalous} for our elliptic curve $A$ (in the sense of Mazur \cite{mazur};  cf. also \cite{local}, Definition 3.2). On the other hand, if one attempts to apply the Serre-Tate expansion method
 in \cite{local} to the case of ${\mathbb G}_m$ instead of $A$ one immediately discovers  that
   ``all primes behave as if they were anomalous;" hence the method is not applicable.
 \end{remark}
 
 \medskip
 
 We next consider the proof of Theorem \ref{thm2}. We need to review some notation from \cite{euler} first. We may and will assume $a_1\neq 0$.
 Consider  two more indeterminates $x,y$,  and consider the  polynomial 
   $$F=F(z_1,z_2,x)\in {\mathbb C}[z_1,z_2][x]$$
   defined by 
      \begin{equation}
     \label{quartic in one variable}
     F:=((a_2-a_3)x^2+z_1-a_2z_2)((a_3-a_1)x^2-z_1+a_1z_2).
     \end{equation}
   For any $c=(c_1,c_2)\in {\mathbb C}^2$ set
   \begin{equation}
   E'_c:=\operatorname{Spec}\ {\mathbb C}[x,y]/(y^2-F(c_1,c_2,x)).
   \end{equation}
   Then we have a morphism
     \begin{equation}
     \label{isogeny}
    \pi: E_c\ra E'_c,\end{equation}
     given by
     $$ x\mapsto x_3,\ \ \ y\mapsto (a_1-a_2)x_1x_2.$$
         Note that, under the assumption that $N(c_1,c_2)\neq 0$, the
         discriminant  of $F$ is in ${\mathbb C}^{\times}$ so  $E'_c$ is a
         smooth plane curve whose smooth projective model ${\mathcal E}'_c$ is
         an elliptic curve and hence \ref{isogeny} is induced by a degree-two
         isogeny of elliptic curves,  
     \begin{equation}
     \label{completed isogeny}
     {\mathcal E}_c\ra {\mathcal E}'_c.\end{equation}
    We also need some well-known formulas relating quartic to cubic
    equations. We recall that  
    the complex plane cubic
    $$y^2=Ax^4+Cx^2+E$$
    is birational  to the Weierstrass equation
    $$v^2=u^3+A_2u^2+A_4u+A_6$$
    where
    $$
    \begin{array}{rcl}
    A_2 & = & C,\\
    \  & \  & \  \\
    A_4 & = & -  4AE,\\
    \  & \ & \  \\
    A_6 & =  & -  4ACE;
    \end{array}
    $$
    the coordinates are related by
    $$x=\frac{2\sqrt{E}(u+C)}{v},\ \ \ y=-\sqrt{E}+\frac{x^2u}{2\sqrt{E}}.$$
    On the other hand recall that the discriminant $\Delta$ of the above elliptic curves
    is given by
    \begin{equation}
    \label{jjj}
    \begin{array}{rcl}
    1728 \Delta & = & (16A_2^2-48A_4)^3-(64A_2^3-288A_2A_4+864A_6)^2\\
    \  & \ & \  \\
    \  & = & (16C^2+192AE)^3-(64C^3-2304ACE)^2,\end{array}
    \end{equation}
    and the $j$-invariant is given by
    \begin{equation}
    \label{jjjj}
    j=\frac{(16C^2+192AE)^3}{\Delta}.
    \end{equation}
So the $j$-invariant of ${\mathcal E}'_c$ is given by \ref{jjjj} with $\Delta$ as in \ref{jjj} and 
\begin{equation}
\label{ACE}
\begin{array}{rcl}
A & = & (a_2-a_3)(a_3-a_1),\\
\  & \  & \ \\
C & = & (a_2-a_3)(a_1c_2-c_1)+(a_3-a_1)(c_1-a_2c_2),\\
\  & \ & \ \\
E & = & (c_1-a_2c_2)(a_1c_2-c_1).
\end{array}
\end{equation}

 \medskip
 
 {\it Proof of Theorem \ref{thm2}}. 
 Recall that we assumed $a_1\neq 0$.
 In view of Theorem \ref{thm1} it is enough to show that, in the notation of Theorem \ref{thm2} and with notation as in the discussion preceding this proof, the following holds: for any root of unity $c_2 \in {\mathbb C}$ 
 we have
 \begin{equation}
 j({\mathcal E}_{(c_1,c_2)})\ra \infty\ \ \ \text{as}\ \ \ c_1 \ra a_1c_2.
 \end{equation}
 
 Indeed if this is the case, with the  root of unity $c_2\in {\mathbb C}$ fixed,  we can apply Theorem \ref{thm1} to the elliptic fibration 
 obtained by compactifying the family of smooth elliptic curves,
 $$({\mathcal E}_{(t,c_2)})_t,\ \ \ t\in {\mathbb C}\backslash \{a_1c_2, a_2c_2, a_3c_2\};$$
 we get that if $p>>0$ then there are only finitely many pairs $(c_1,c_2)$, with $c_1\in T_p$, such that 
 ${\mathcal E}_{(c_1,c_2)}$ is CL. 
 
 To conclude the proof
 fix a fundamental domain $F$ for $SL_2(\bZ)$ in the complex upper half plane 
 that contains all complex numbers with real part in $(-1,1)$ and imaginary part in $(1,\infty)$, say.
 If $\tau\in F$ corresponds to ${\mathcal E}'_c$ then  ${\mathcal E}_c$ (which admits an isogeny of degree $2$ to ${\mathcal E}'_c$) corresponds 
 to one of the numbers $2\tau,2\tau+1,\tau/2$. So it is enough to prove that
  \begin{equation}
 j({\mathcal E}'_{(c_1,c_2)})\ra \infty\ \ \ \text{as}\ \ \ c_1 \ra a_1c_2.
 \end{equation}
 Now, using \ref{jjj}, \ref{jjjj}, \ref{ACE}, one gets  that, as $c_1 \ra a_1c_2$, we have
 $$E\ra 0,\ \ 
 C\ra (a_3-a_1)(a_1-a_2)c_2\neq 0,\ \ \ \ \Delta\ra 0,$$
so $j\ra \infty$ and we are done.
 \qed

  \bigskip
  
  \begin{remark}
    Theorem \ref{thm2} states that for any choice of
   $(a_1,a_2,a_3)$ there exists a bound such that for any prime larger than that bound one can find ``many" pairs $(c_1,c_2)$ such that a certain property holds. One can turn the tables and show that there exists a  pair $(c_1,c_2)$ for which 
   one can find ``many" triples $(a_1,a_2,a_3)$ and ``many" primes such that that same  property holds. The latter is easier and does not require 
  an input from Theorem \ref{thm1}. Here is an example of such a statement.
  
  \medskip
   
   Let $m\in {\mathbb Z}$ and set
   $a_1=2,\ \ a_2=0,\ \ a_3=m.$

   \medskip
   
   {\it Claim}. For all but finitely many $m$ 
   the following holds: there are infinitely many primes $p$ such that 
   
   1) $N(1,1)\not\in R_p^{\times}$;
   
   2) ${\mathcal E}_{(1,1)}$ has ordinary reduction but is  not CL.
   
   \medskip
   
   To check the Claim note that
   $$F(1,1,x)=m(2-m)x^4-2x^2+1,$$
   so the affine elliptic curve defined by
   \begin{equation}
   \label{boba1}
   y^2=F(1,1,x)\end{equation}
   is birationally equivalent over ${\mathbb C}$ to the affine elliptic 
   curve defined by
   \begin{equation}
   \label{boba2}
   y^2=x^3-2x^2-4m(2-m)x+8m(2-m).\end{equation}
   The latter has $j$-invariant given by
   $$j(m)=\frac{P_6(m)}{P_6(m)-P_4(m)}$$
   where $P_6,P_4\in {\mathbb Z}[t]$ are polynomials of degree $6$ and $4$  respectively. Since there are only finitely many numbers   in ${\mathbb Q}$ that appear as  $j$-invariants of
    CM elliptic curves there is a cofinite set of integers $m$ such 
   the curve defined by \ref{boba2}
   is not a CM curve. 
   Fix such an $m$.
   Then the curve defined by \ref{boba1}, and hence the curve ${\mathcal E}_c$,  is not a CM curve.
   Now there are infinitely many primes $p$ such that ${\mathcal E}_c$ has ordinary reduction, i.e., for $A_{p-1}$ the Hasse invariant as in \cite{euler}, 
   $$A_{p-1}(1,1)\in {\mathbb Z}_{p}^{\times}.$$
   Finally there are infinitely many primes $p$ such that, in addition, 
   $$N(1,1)=\prod_{i=1}^3(1-a_i)=m-1\in {\mathbb Z}_{p}^{\times}.$$ 
   For such primes ${\mathcal E}_c$ is not a CL curve.
   This ends the proof of our Claim.
   \end{remark}


\begin{thebibliography}{AA}
  
    
  
  
\bibitem{neron} S. Bosch, W. L\"{u}tkebohmert, M. Raynaud, {\it N\'{e}ron
  Models}, Springer, 1980. 

\bibitem{char} A. Buium, {\it Differential characters of
abelian varieties over $p$-adic fields},  Invent. Math. 122, 2
(1995), 309-340.

\bibitem{euler} A. Buium, E. Previato,
{\it Arithmetic Euler top},   J. Number Theory, 173 (2017),  37-63.


\bibitem{local} A. Buium, B. Poonen, {\it  Independence of points on elliptic
curves arising from special
 points on modular and Shimura curves, II: local results},
 Compositio Math., 145 (2009), 566-602.

\bibitem{bogo} K. Kodaira, {\it On compact analytic surfaces II, III}, Annals of Math. (2) 77 (1963), 563-626, and 78 (1963), 1-40.

\bibitem{mazur} Mazur, B., {\it Rational points of abelian
varieties with values in towers of number fields}, Invent.\ Math. 18
(1972), 183-266.

\end{thebibliography}
\end{document}